\documentstyle[12pt,amsfonts]{article}

\newcommand{\E}{{\Bbb E}}

\newcommand{\h}{{\cal H}}

\newcommand{\R}{\Bbb R}

\begin{document}
\title{Integration by parts formulae for degenerate diffusion measures on 
path spaces and diffeomorphism groups}

\author{K.D.  Elworthy \and Yves  Le Jan \and Xue-Mei Li }

\date{}

\maketitle

\footnote{Research supported by EPSRC grant GR/H67263 and EC programme
Science plan ERB 4002PL910459.}

\begin{abstract}
Integration by parts formulae are given for a class of measures on the space of
paths of a smooth manifold $M$ determined by the laws of degenerate diffusions.
The mother of such formulae, on the path space of diffeomorphism group of $M$
is shown to arise from a quasi-invariance property of measures determined
by stochastic flows. From this the other formulae are derived by
 filtering out redundant noise using an associated LeJan-Watanabe connection.
\end{abstract}

Resum\'e: Des formules d'int\'egration par parties sont donn\'ees pour
les mesures sur l'espace des chemins associe\'es  des E.D.S,
eventuellement d\'eg\'en\'er\'ees. La formule m\'ere est donn\'ee
sur l'espace des chemins dans le groupe des diff\'eomorphismes
de la vari\'et\'e. Elle proc\'ede de la quasi invariance des
lois du flot stochastique. La formule sur l'espace des chemins
qui en d\'erive est obtenue en filtrant le bruit redondant 
\`a l'aide de la connexion intraduite par S. Watanabe et le deuxi\`eme
auteur.

\bigskip

{\bf A.} Consider the equation
\begin{equation}\label{1}
dx_t=X(x_t)\circ dB_t+A(x_t)dt
\end{equation}
on a compact $C^\infty$ manifold $M$. Here $\{B_t: t\ge 0\}$ is a
Brownian motion on $\R^m$ and $A$ is a $C^\infty$ vector field
while $X: M\times \R^m\to TM$ is $C^\infty$ and gives a linear map
$X(x,\cdot)=X(x)(\cdot): \R^m\to T_xM$ for each $x\in M$.

\bigskip

Let $\{\xi_t: t\ge 0\}$ be the solution flow to (\ref{1}) with
 $T_x\xi_t: T_xM \to T_{\xi_t(x)}M$ its derivative map at $x\in M$. For fixed 
$T>0$ suppose $k$ is an element of the Cameron-Martin space
 $H=L_0^{2,1}([0,T];\R^m)$ and consider the random, time dependent O.D.E.
on $M$, parametrized by $\tau \in \R$:
\begin{equation} \label{2}
\begin{array}{llll}
{\partial \over \partial t}H_t^\tau(x)&=&\tau \left(T_x\xi_t\right)^{-1}
X\left(\xi_t(H_t^\tau(x))\right) \dot k_t,  \hskip 12pt & 0\le t\le T,\\
H_0^\tau(x)&=&x, &x\in M.
\end{array}
\end{equation}

The solution exists for $0\le t \le T$ and we perturb our flow to $\xi_t^\tau$
given by 
$$\xi_t^\tau(x)=\xi_t(H_t^\tau(x)), \hskip 30pt 0\le t\le T.$$
Then $\{\xi_t^\tau(x): 0\le t\le T\}$ is a semi-martingale satisfying

\begin{equation}\label{3}
d\xi_t^\tau(x)=X\left(\xi_t^\tau(x)\right)\circ dB_t
+A(\xi_t^\tau(x))dt+\tau X(\xi_t^\tau(x))\dot k_tdt
\end{equation}

\noindent
(It is given by replacing $B_\cdot$ by $B_\cdot+\tau k$ in (\ref{1}))
and so the laws of $\xi_\cdot^\tau(x)$, $\tau \in \R$, are mutually 
equivalent by the Cameron-Martin theorem. By differentiating the
resulting two expressions for $\E F(\xi_\cdot^\tau(x))$ with respect
to $\tau$ at $\tau=0$ in the usual way, e.g. as in \cite{EL-LIibp} we see

\bigskip

{\it Let $F: C_x\left([0,T];M\right) \to \R$ be $C^1$ with bounded
derivative. Then}
\begin{equation}\label{4}
\E dF\left(T\xi_\cdot \int_0^\cdot \left(T\xi_s\right)^{-1}
\left(X(\xi_s(x)\dot k_s \right)ds   \right)
=  \E F(\xi_\cdot (x)) \int_0^T <\dot k_s, dB_s>_{\R^m}.
\end{equation}

This formula could also be deduced from the usual integration by parts
formula on $C_0\left([0,T];\R^m\right)$ as in \cite{Bismut81b}. However
the above argument also shows that
 {\it the law $\mu_{\cal D}$ of
$\xi_\cdot$ on the space of paths $C_{id} \left ([0, T]; \hbox{Diff}M \right )$ 
into the group of diffeomorphisms of $M$ is quasi-invariant under
 the transformations $\xi_\cdot \to \xi_\cdot(H_\cdot^\tau)$ with
$H_\cdot^\tau$ as above},  $\tau\in \R$, $k_\cdot\in H$.
 Let ${\cal H}$ be the
reproducing kernel Hilbert space of vector fields on $M$ for the laws
of the random field $X(\cdot)B_1$. Then $K_\cdot(-)=X(-)k_\cdot$
is an arbitrary element of $L^{2,1}_0 \left([0,T];{\cal H}\right)$.
Let $\{W_t: 0\le t \le T\}$ be the martingale part of 
$\{\int_0^t d\xi_s\circ \xi_s^{-1}: 0\le t\le T\}$. Then
$dW_t=X(\cdot)dB_t$
and we obtain:

\begin{equation}\label{5}
\begin{array}{l}
\int_{C_{id}\left([0,T]; \hbox{Diff}M\right)}
dF\left(T\theta_\cdot \int_0^\cdot ad(\theta_s^{-1}) \dot K_s ds
\right)\mu_{\cal D}(d\theta)\\
=\int_{C_{id}\left([0,T];\hbox{Diff}M\right)}
F(\theta_\cdot) \int_0^T <\dot K_s, dW_s>
_{\cal H}  \mu_{\cal D}(d\theta)\\
\end{array}
\end{equation}
{\it for any $C^1$ map $F: C_{id}\left([0,T]; \hbox{Diff}M\right) \to \R$
with suitably bounded derivative } $dF$. Here $ad(\theta_s^{-1})$
refers to the adjoint action of $\theta_s^{-1}$ on the space of vector
fields of $M$ considered as the Lie algebra of Diff M.

\bigskip

\noindent
{\bf Remarks:}  (i): In (\ref{4}), (\ref{5}) and below we could
equally
well take $k_\cdot$ or $K_\cdot$ to be the adapted processes with
sample paths in $L_0^{2,1}$, and in $L^{1+\epsilon}$ as
$L_0^{2,1}$-valued functions for some $\epsilon>0$, c.f \cite{EL-LIibp}.

(ii). for a general stochastic flow of diffeomorphisms of adequate
regularity formula (\ref{5}) remains valid but ${\cal H}$ may be infinite
dimensional. Indeed there is an analogous formula for any right
invariant system on a Hilbert manifold with sufficiently regular group
structure.

(iii). Equation (\ref{2}) defining $H_t^\tau, 0\le t\le T$, is just
the right invariant equation on Diff$M$  corresponding to the element
$ad(T\xi_t)^{-1}(\dot K_t)$. In general it is a function of
$\xi_\cdot$, but if there exists $V$ in ${\cal H}$ which commutes
with all other elements of ${\cal H}$ then setting $K_t=tv$ we obtain
an $H_\cdot^\tau$, right multiplication by which leaves $\mu_{\cal D}$
quasi-invariant .

\bigskip

{\bf B.} Equation(\ref{4}) is not an integration by parts formula on
$C_x([0,T];M)$, strictly speaking, since we do not have a vector field
on $C_x([0,T];M)$ and there is possibly redundant noise coming from
lack of injectivity of $X$. To put (\ref{4})  in correct form we
need the conditional expectation

\begin{equation}\label{6}
\bar V_t^k;=\E\left\{ \left.   T\xi_t\int _0^t
(T\xi_s)^{-1}\left(X(\xi_s(x)\dot K_s\right)ds\,\right|\,\sigma\left\{\xi_s(x):
0\le s \le T\right\}   \right\}.
\end{equation}

For simplicity here we assume that $X(x)$ has constant rank
(non-singularity). Then $X$ has image a subbundle $I(X)$ of $TM$ on
which $X$ induces (i) a Riemannian metric $\{<,>_x, x \in M\}$ and
(ii) an affine connection $\check \nabla $, which is a metric
connection for this metric (the LeJan-Watanabe connection). If 
$Y(x): I(X)_x\to \R^m$
is the adjoint of $X(x):\R^m\to I(X)_x$ these satisfy
$$X(x)Y(x)v=v, \hskip 14pt v\in T_xM$$
and $$\check \nabla_v Z=X(x) d\left(Y(\cdot)Z(\cdot)\right)v$$
for all $v\in T_xM$ and $Z$ a section of $I(X)$. There is also an
'adjoint semi-connection' $\hat \nabla$ to $\check \nabla$ which gives
a derivative $\hat \nabla_{Z_2}Z_1$ of a vector field $Z_1$ 
 in the direction of a section $Z_2$ of $I(X)$ by
$$\hat \nabla_{Z_2}Z_1=\check \nabla_{Z_1}Z_2 -[Z_1,Z_2].$$
Let $v_0\in T_xM$ and set 
$$\bar v_t=\left\{ \left.  T\xi_t(v_0) \right| \sigma\left\{\xi_s(x):
0\le s \le T\right\}   \right\}.$$
Then, from \cite{EL-LJ-LI2} following \cite{EL-LJ-LI} and
\cite{EL-YOR}, if $A(x)\in I(X)_x$ for each $x\in M$ then 
$\{\bar v_t: 0\le t\le T\}$ satisfies

\begin{equation}\label{7}
{\hat D \bar v_t\over \partial t} =-{1\over 2}
\check{\hbox{Ric}}^{\#}(\bar v_t)
 + \check\nabla A(\bar v_t)
\end{equation}
where $<\check{\hbox{Ric}}^{\#}(u_1), u_2> 
=\sum_{i=1}^r <\check R(e_i, u_1)u_2, e_i>_x$ for
$$\check R: TM\bigoplus TM \to {\cal L} \left(I(X); I(X)\right)$$
the curvature tensor of $\check \nabla $ and $e_1, \dots, e_r$ an
orthonormal basis for $I(X)_x$ where $u_1\in T_xM$, $u_2\in I(X)_x$
and $I(X)_x$ is identified with its image in $T_xM$. This can also be
written as

\begin{equation}\label{8}
D\bar v_t=\nabla X(\bar v_t)\left(\tilde {//_t}\circ d \tilde
B_t\right)+
\nabla A(\bar v_t)dt-{1\over 2} \check{\hbox{Ric}} (\bar v_t)dt
\end{equation}
using the Levi-civita connection $\nabla$ for any Riemannian metric on
$M$ which extends that of $I(X)$. Here 
$\tilde {//_t}\circ d{\tilde B_t}
=Y\left(\xi_t(x)\right)  X\left(\xi_t(x)\right)\circ dB_t$ and can
be represented by a parallel translation of the differential d$\check
B_t$ of the martingale part of the stochastic anti-development of
$\{\xi_t(x): 0\le t\le T\}$. Finally we have, from (\ref{4}):
{\it        For any given $v_0\in T_xM$ let $\{W_t^A(v_0): 0\le t\le T\}$ be the
solution of (\ref{8}), or equivalently of (\ref{7}) if $A(y)\in I(y)$
for each $y\in M$. Let $\mu_x$ be the law of
 $\{\xi_t(x): 0\le t \le T\}$ on $C_x\left([0,T];M\right)$. Then 

\begin{equation}\label{9}
\begin{array}{l}
\int_{C_x\left([0,T];M \right)} dF\left( W_\cdot^A 
\int_0^\cdot \left(W_s^A\right)^{-1}
\left(X(\sigma(s)) \dot K_s\right) ds \right) \mu_x(d\sigma))\\
=\int_{C_x\left([0,T];M\right)} F(\sigma)\int_0^T
\left<X(\sigma(s))\dot K_s, \check{//_s} d\check B_s
\right>_{T_{\sigma(s)}M}  \mu_x(dx)
\end{array}
\end{equation}
for any $C^1$ map $F: C_x\left([0,T];M\right)\to \R$ with bounded derivative
and any $K_\cdot \in L_0^{2,1}([0,T];\R^m)$.}

\bigskip

{\bf Remarks:} (i) For the law of $\xi_\cdot(x)$ for a more general
stochastic flow of diffeomorphisms, as described in \S A formula
(\ref{9}) remains valid with $X(\sigma(s))\dot k_s$ replaced by
 $\dot K_s(\sigma(s))$.

(ii). $k$ (or $K$) can be non-anticipating processes as in Remark (i) 
of \S A.

(iii). For an analogous approach for non-degenerate equations see
\cite{EL-LIibp}, for LeJan-Watanabe connections in the non-degenerate
case see \cite{EL-LJ-LI}; the degenerate case with detailed
discussions of these results and of the singular case will appear in
\cite{EL-LJ-LI2}. 

(iv) Following on from remark (ii) in section \S A note this
formulation
extends to loop spaces if one adapts the point of view used by Driver
in \cite{Driver95} for loops on Lie groups. Replacing $B_t$ by a
centred Gaussian process $B_{t,\alpha}, t\in \R^+, \alpha \in (0,1)$
with covariance $\E\left(B_{t,\alpha}^i B_{s, \beta}^j\right)=\delta_{i.j}
s\wedge t(\alpha\wedge \beta -\alpha \beta)$, one gets a process
$\xi_{t,\alpha}$ which represents a Brownian motion in the group of
based loops in the diffeomorphism group, acting naturally on pinned
loops of $M$. Taking into account the $\alpha$ dependence, integration
by parts formulae (\ref{4}), (\ref{5}) and (\ref{9}) extend naturally
to this situation.

\vskip 30pt

\begin{center}{\large  Version francaise abre\'ege\'e}
\end{center}

{A.} Considerons l'equation diff\`erentielle stochastique sur une
vari\'et\'e compacte $M$. $\{B\}$ designe un mouvement brownien
sur $\R^m$ et $A$ et un champ de vecteurs $C^\infty$ tandis que 
$X: \R^m\times M \to TM$ est $C^\infty$ et donne une application
lineaire $X(x): \R^m \to T_xM$ pour tout $x\in M$. 

Soit $\{\xi_t^\tau\}$ le flot stochastique solution de

\begin{equation}\label{1f}
dx_t=X(x_t)\circ dB_t+A(x_t)dt
\end{equation}
 et $T_x\xi_t: T_xM \to T_{\xi_t(x)}M$ sa deriv\'ee en $x\in M$.
 Pour $T>0$ fix\'e, soit $k$ un 
element de l'espace de Cameron Martin $H=L_0^{2,1}([0,T];\R^m)$.
Considerons l'equation diff\'erentielle ordinaire sure $M$
param\'etr\'ee par $\tau \in \R$.
\begin{equation} \label{2f}
\begin{array}{llll}
{\partial \over \partial t}H_t^\tau(x)&=&\tau \left(T_x\xi_t\right)^{-1}
X\left(\xi_t(H_t^\tau(x))\right) \dot k_t,  \hskip 12pt & 0\le t\le T,\\
H_0^\tau(x)&=&x, &x\in M.
\end{array}
\end{equation}
\bigskip

La solution existe pour $0\le t\le T$. Definissons le flot perturb\'e
$\xi_t^\tau$ par
 $$\xi_t^\tau(x)=\xi_t(H_t^\tau(x)), \hskip 30pt 0\le t\le T.$$
Alors $\{\xi_t^\tau\}$ est une semi-martingale v\,erifiant (3). 
Les lois des $\xi_\cdot^\tau(x)$ sont equivalentes. En
diff\'erentiant par rapport \`a $\tau$ en $\tau=0$ (cf par
exemple \cite{EL-LIibp}) on voit que:

{\it Soit
 $F: C_x\left([0,T];M\right) $ suppos\'ee $C^1$ et \`a deriv\'ees
boune\'es. Alors}
\begin{equation}\label{4f}
\E dF\left(T\xi_\cdot \int_0^\cdot \left(T\xi_s\right)^{-1}
\left(X(\xi_s(x)\dot k_s \right)ds   \right)
=  \E F(\xi_\cdot (x)) \int_0^T <\dot k_s, dB_s>_{\R^m}.
\end{equation}

Cette formule peut \^etre deduite de la formule usuelle d'integration
par parties sur  $C_0\left([0,T];\R^m\right)$ comme dans \cite{Bismut81b}. 
Cependant, l'argument pr\'ec\'edent montre aussi que {\it  la loi
$\mu_{\cal D}$ de $\xi_\cdot$ sur l'espace des trajectoines 
 $C_{id} \left ([0, T]; \hbox{Diff}M \right )$  dans le groupe
des diffeomorphismes de $M$ est quasi invariante par les
transformations   $\xi_\cdot \to \xi_\cdot(H_\cdot^\tau)$, o\`u
$H_\cdot$ est d\'efini ci dessus,  $\tau\in \R$, $k_\cdot\in H$.}
Soit ${\cal H}$ l'espace auto--reproduisant de champs de
vecteurs sur $M$ associ\'e au champ gaussien 
 $X(\cdot)B_1$. $K_\cdot(-)=X(-)k_\cdot$ est alors un element
arbitraire de $L^{2,1}_0 \left([0,T];{\cal H}\right)$.
Soit $\{W_t: 0\le t \le T\}$ le partie martingale de
$\{\int_0^t d\xi_s\circ \xi_s^{-1}: 0\le t\le T\}$. 
Alores   $dW_t=X(\cdot)dB_t$ et nous obtenons:

\begin{equation}\label{5f}
\begin{array}{l}
\int_{C_{id}\left([0,T]; \hbox{Diff}M\right)}
dF\left(T\theta_\cdot \int_0^\cdot ad(\theta_s^{-1}) \dot K_s ds
\right)\mu_{\cal D}(d\theta)\\
=\int_{C_{id}\left([0,T];\hbox{Diff}M\right)}
F(\theta_\cdot) \int_0^T <\dot K_s, dW_s>
_{\cal H}  \mu_{\cal D}(d\theta)\\
\end{array}
\end{equation}
{\it pour toute application $C^1$  $F: C_{id}\left([0,T]; \hbox{Diff}M\right) \to \R$
poss\'edant une deriv\'ee ad\'equatement born\'ee.}
On d\'esigne par $ad(\theta_s^{-1})$ la repr\'esentation adjointe
de f $\theta_s^{-1}$ sur l'espace des champs de vecteurs sur 
 $M$ consid\'er\'e comme l'alg\'ebre de Lie  de Diff M.

\bigskip

\noindent
{\bf B.} L'\'equation (\ref{4f}) n'est pas une formule d'int\'egration
par parties sur $C_x([0,T];M)$, \`a proprement parler. Pour mettre
(\ref{4f}) sous une forme plus appraprie\'e, il faut prendre
l'esp\'erance conditionelle

\begin{equation}\label{6f}
\bar V_t^k;=\E\left\{ \left.   T\xi_t\int _0^t
(T\xi_s)^{-1}\left(X(\xi_s(x)\dot K_s\right)ds\,\right|\,\sigma\left\{\xi_s(x):
0\le s \le T\right\}   \right\}.
\end{equation}

\bigskip

Pour simplifier, nous supposerons que le rang de $X(x)$ est constant.
Alors, l'image de $X$ est un sous-fibr\'e $I(X)$ de $TM$ sur lequel
$X$ induit (i) une m\'etrique riemanniene   $\{<,>_x, x \in M\}$    et 
(ii) une connexion affine   $\{\check \nabla$, qui preserve
la m\'etrique en question. Si nous notons 
$Y(x): I(X)_x\to \R^m$
l'application adjointe de  $X(x):\R^m\to I(X)_x$, on a 
$$X(x)Y(x)v=v, \hskip 14pt v\in T_xM$$
et $$\check \nabla_v Z=X(x) d\left(Y(\cdot)Z(\cdot)\right)v$$
pout tout  $v\in T_xM$ et pour  toute section  $Z$ de $I(X)$. 
Il y a aussi une 'semi connection adjointe' $\hat \nabla$
de $\breve \nabla$ qui associe une derivee
$\hat \nabla_{Z_2}Z_1$, \`a tout champ de vecteur   $Z_1$ 
dans  la direction $Z_2$ d'une section de $I(X)$, 
d\'efinie par
$$\hat \nabla_{Z_2}Z_1=\check \nabla_{Z_1}Z_2 -[Z_1,Z_2].$$
Si  $v_0\in T_xM$ posons
$$\bar v_t=\left\{ \left.  T\xi_t(v_0) \right| \sigma\left\{\xi_s(x):
0\le s \le T\right\}   \right\}.$$
Alors, d'apr\'es  \cite{EL-LJ-LI2} et suite \`a  \cite{EL-LJ-LI} et
\cite{EL-YOR}, si $A(x)\in I(X)_x$  pour tout  $x\in M$,
$\{\bar v_t: 0\le t\le T\}$ verifie

\begin{equation}\label{7f}
{\hat D \bar v_t\over \partial t} =-{1\over 2}
\check{\hbox{Ric}}^{\#}(\bar v_t)
 + \check\nabla A(\bar v_t)
\end{equation}
o\`u  $<\check{\hbox{Ric}}^{\#}(u_1), u_2> 
=\sum_{i=1}^r <\check R(e_i, u_1)u_2, e_i>_x$,

$$\check R: TM\bigoplus TM \to {\cal L} \left(I(X); I(X)\right)$$
\'etant le tenseur de courbure de  $\check \nabla $ et $e_1, \dots,
e_r$  une base orthonormale de  $I(X)_x$:  $u_1$ est un element
de $ T_xM$, $u_2$ un element de $I(X)_x$ et  $I(X)_x$  est
identifi\'e \`a son image dans  $T_xM$. On peut aussi l'\'ecrire
sous la forme

\begin{equation}\label{8f}
D\bar v_t=\nabla X(\bar v_t)\left(\tilde {//_t}\circ d \tilde
B_t\right)+
\nabla A(\bar v_t)dt-{1\over 2} \check{\hbox{Ric}} (\bar v_t)dt
\end{equation}
en utilisant la connection de  Levi Civita associ\'ee \`a
toute metrique riemanniene sur $M$ \'etendant celle d\'efinie
sur $I(X)$. $\tilde {//_t}\circ d{\tilde B_t}
=Y\left(\xi_t(x)\right)  X\left(\xi_t(x)\right)\circ dB_t$
peut \^etre represent\'e par transport parall\`ele de la differentielle
 d$\check B_t$ de la partie martingale de l'antideveloppement
stochastique de $\{\xi_t(x): 0\le t\le T\}$. Enfin on deduit de
(\ref{4f}) le r\'esultat suivant:

{it Pour tout   $v_0\in T_xM$ soit  $\{W_t^A(v_0): 0\le t\le T\}$ 
la solution de  (\ref{8f}) (ou de (\ref{7f}) si  $A(y)\in I(y)$
pour tout $y\in M$). Alors

\begin{equation}\label{9f}
\begin{array}{l}
\int_{C_x\left([0,T];M \right)} dF\left( W_\cdot^A 
\int_0^\cdot \left(W_s^A\right)^{-1}
\left(X(\sigma(s)) \dot K_s\right) ds \right) \mu_x(d\sigma))\\
=\int_{C_x\left([0,T];M\right)} F(\sigma)\int_0^T
\left<X(\sigma(s))\dot K_s, \check{//_s} d\check B_s
\right>_{T_{\sigma(s)}M}  \mu_x(dx)
\end{array}
\end{equation}
pour toute application  $C^1$  $F$ de $ C_x\left([0,T];M\right)$
a deriv\'ee born\'ee et tout  $K_\cdot \in L_0^{2,1}([0,T];\R^m)$.}

\noindent {\bf Addresses:}\\

\noindent
 K. D. Elworthy:  Mathematics Department, University of
 Warwick, Coventry CV4 7AL, UK.\\
Yves, Le Jan:  D\'epartment de Math\'ematique, Universit\'e Paris Sud,\\
 91405 Orsay, France\\
Xue-mei Li: Department of Mathematics, U-9,
 University of Connecticut, 196 Auditorium road, Storrs, Connecticut
 06269-3009, USA\\


\bigskip

\begin{thebibliography}{ELJL95}

\bibitem[Bis81]{Bismut81b}
J.~M. Bismut.
\newblock Martingales, the {M}alliavin calculus and {H}\h{o}rmander's theorem.
\newblock In D.~Williams, editor, {\em Stochastic Integrals, {L}ecture {N}otes
  in {M}aths. 851}, pages 85--109. Springer-Verlag, 1981.

\bibitem[Dri95]{Driver95}
Bruce~K. Driver.
\newblock Towards calculus and geometry on path spaces.
\newblock In {\em Stochastic Analysis: AMS Proceedings of symposium in pure
  Math. Series}, pages 423--426. AMS. Providence, Rhode Island, 1995.

\bibitem[EL96]{EL-LIibp}
K.D. Elworthy and Xue-Mei Li.
\newblock A class of integration by parts formulae in stochastic analysis i.
\newblock In {\em {I}t\^o's Stochastic Calculus and Probability Theory
  (dedicated to Prof. It\^o on the occasion of his eightieth birthday)}.
  Springer-Verlag, 1996.

\bibitem[ELJL95]{EL-LJ-LI}
K.~D. Elworthy, Yves Le~Jan, and Xue-Mei Li.
\newblock Concerning the geometry of stochastic differential equations and
  stochastic flows.
\newblock To appear in 'New Trends in stochastic Analysis', Proc. Taniguchi
  Symposium, Sept. 1995, Charingworth, ed. K. D. Elworthy and S. Kusuoka, I.
  Shigekawa, World Scientific Press, 1995.

\bibitem[ELL]{EL-LJ-LI2}
K.D. Elworthy, Yves LeJan, and Xue-Mei Li.
\newblock In preparation.

\bibitem[EY93]{EL-YOR}
K.~D. Elworthy and M.~Yor.
\newblock Conditional expectations for derivatives of certain stochastic flows.
\newblock In J.~Az\'ema, P.A. Meyer, and M.~Yor, editors, {\em Sem. de Prob.
  XXVII. Lecture Notes in Maths. 1557}, pages 159--172. Springer-Verlag, 1993.

\end{thebibliography}
\end{document}